\documentclass[11pt]{article}
\usepackage{amscd,amssymb,epic,amsmath}
\usepackage[all]{xy}

\newtheorem{proposition}{Proposition}

\newtheorem{lemma}{Lemma}
\newtheorem{corollary}{Corollary}
\newtheorem{remark}{Remark}

\title{\Large\bf Spherical principal series of quantum Harish-Chandra
modules}\date{}

\begin{document}
\maketitle

\centerline{O.Bershtein$^\dagger$, A.Stolin$^\ddagger$,
L.Vaksman$^\dagger$} \centerline{$^\dagger$Institute for Low
Temperature Physics and Engineering, Kharkov, Ukraine}
\centerline{$^\ddagger$Chalmers University of Technology,
G\"{o}teborg, Sweden} \centerline{bershtein@ilt.kharkov.ua,
astolin@math.chalmers.se, vaksman@ilt.kharkov.ua}

\bigskip

\begin{abstract}
The non-degenerate spherical principal series of quantum
Harish-Chandra modules is constructed. These modules appear in the
theory of quantum bounded symmertic domains.
\end{abstract}



\section{Introduction}

In \cite{Helg1} the unit disc $\mathbb D=\{z \in \mathbb C| |z|<1\}$
is considered as the Poincare model of hyperbolic plane. The
Plancherel formula for $\mathbb D$ is one of the most profound
results of non-commutative harmonic analysis, and $q$-analogs of this
result are well-known for $q \in (0,1)$ being the deformation
parameter, see, e.g., \cite{VakShkl}. It is important to note that
representations of spherical principal series is a crucial tool in
decomposing quasi-regular representation both in the classical and
quantum case.

The unit disc is the simplest bounded symmetric domain. The
Plancherel formula is known for any bounded symmetric domain in the
classical setting. On the other hand, it is absent in the quantum
case. One of the obstacles here is some difficulties in producing
non-degenerate spherical principal series of Harish-Chandra modules
over a quantum universal enveloping algebra. In this paper we
overcome these difficulties. A geometrical approach to the
representation theory is used instead of the traditional construction
of principal series, which is inapplicable in the quantum case.
Hence, we generalize the results of \cite{SStV}.

The Casselman's theorem claims that any simple Harish-Chandra module
can be embedded in a module of non-degenerate principal series. Thus,
one has another class of applications for the constructions of this
paper, beyond the harmonic analysis.

\section{A quantum analog of the open \boldmath $K$-orbit in
$B\backslash G$}

Let $(a_{ij})_{i,j=1,\ldots,l}$ be a Cartan matrix of positive type,
$\mathfrak g$ the corresponding simple complex Lie algebra. So the
Lie algebra can be defined by the generators $e_i,f_i,h_i,i=1,...,l$,
and the well-known relations, see \cite{Jantzen}. Let $\mathfrak{h}$
be the linear span of $h_i,i=1,...,l$. The simple roots $\{\alpha_i
\in \mathfrak{h}^*|i=1,...,l\}$ are given by $\alpha_i(h_j)=a_{ji}$.
Also, let $\{\varpi_i| i=1,...,l\}$ be the fundamental weights, hence
$P=\bigoplus_{i=1}^l
\mathbb{Z}\varpi_i=\{\lambda=(\lambda_1,...,\lambda_l)|\,\lambda_j
\in \mathbb{Z}\}$ is the weight lattice and $P_+=\bigoplus_{i=1}^l
\mathbb{Z_+}\varpi_i=\{\lambda=(\lambda_1,...,\lambda_l)|\,\lambda_j
\in \mathbb{Z}_+\}$ is the set of integral dominant weights.

Fix $l_0 \in \{1,...,l\}$, together with the Lie subalgebra
$\mathfrak{k} \subset \mathfrak{g}$ generated by
$$
e_i, f_i,\;\; i\neq l_0; \qquad h_i, \;\;i=1,...,l.
$$
Define $h_{0} \in \mathfrak{h}$ by
$$
\alpha_i(h_0)=0, \; i \neq l_0; \qquad \alpha_{l_0}(h_0)=2.
$$
We restrict ourself by Lie algebras $\mathfrak{g}$ that can be
equipped with a $\mathbb Z$-grading as follows:
\begin{equation}\label{par_type}
\mathfrak{g}=\mathfrak{g}_{-1}\oplus\mathfrak{g}_0\oplus\mathfrak{g}_{+1},
\quad \mathfrak{g}_j=\{\xi\in\mathfrak{g}|\:[h_0,\xi]=2j\xi\}.
\end{equation}
Let $\delta$ be the maximal root, and $\delta=\sum_{i=1}^l
c_i\alpha_i$. (\ref{par_type}) holds if and only if $c_{l_0}=1$. In
this case $\mathfrak{g}_0=\mathfrak{k}$ and the pair
$(\mathfrak{g},\mathfrak{k})$ is called a Hermitian symmetric pair.

%Let $W$ be the Weyl group of the root system $R$, $w_0 \in W$ the
%longest element. The irreducible bounded symmetric domain $\mathbb D$
%related to the pair $(\mathfrak{g},\mathfrak{k})$ is a domain of tube
%type if and only if $\varpi_{l_0}=-w_0\varpi_{l_0}$.

Fix a Hermitian symmetric pair $(\mathfrak{g},\mathfrak{k})$. Let $G$
be a complex algebraic affine group with ${\rm Lie}\, (G)
=\mathfrak{g}$ and $K \subset G$ the connected subgroup with ${\rm
Lie}\,(K) =\mathfrak{k}$. Consider the Lie subalgebra $\mathfrak{b}
\subset \mathfrak{g}$ generated by $e_i, h_i, i =1,...,l$, together
with the corresponding connected subgroup $B \subset G$. The
homogeneous space $G/P$ has a unique open $K$-orbit $\Omega \subset
X$. It is an affine algebraic variety which is a crucial tool in
producing non-degenerate principal series of Harish-Chandra modules
(see \cite{Schmid}).

We introduce a $q$-analog of the algebra of regular functions on the
open orbit.

Recall some background and introduce the notation. First of all,
recall some notions from the quantum group theory \cite{Jantzen}. In
the sequel the ground field is $\mathbb C$, $q \in (0,1)$, and all
the algebras are associative and unital.

Denote by $d_i, i=1,\ldots,l$, such positive coprime integers that
the matrix $(d_ia_{ij})_{i,j=1,\ldots,l}$ is symmetric. Recall that
the quantum universal enveloping algebra $U_q \mathfrak{g}$ is a Hopf
algebra defined by the generators $K_i$, $K_i^{-1}$, $E_i$, $F_i$,
$i=1,\ldots,l$ and the relations:
$$K_iK_j=K_jK_i,\quad K_iK_i^{-1}=K_i^{-1}K_i=1,$$
\begin{equation*}
K_iE_j=q_i^{a_{ij}}E_jK_i,\quad K_iF_j=q_i^{-a_{ij}}F_jK_i,
\end{equation*}
$$
E_iF_j-F_jE_i=\delta_{ij}\,\frac{K_i-K_i^{-1}}{q_i-q_i^{-1}},
$$
\begin{equation*}
  \sum\limits_{m=0}^{1-a_{ij}}(-1)^m\left[1-a_{ij} \atop m\right]_{q_i}
  E_i^{1-a_{ij}-m}E_jE_i^m=0,
\end{equation*}
\begin{equation*}
  \sum\limits_{m=0}^{1-a_{ij}}(-1)^m\left[1-a_{ij} \atop m\right]_{q_i}
  F_i^{1-a_{ij}-m}F_jF_i^m=0,
\end{equation*}
where $q_i=q^{d_i}$, $1\leq i\leq l$, and
\begin{equation*}
  \left[m \atop n\right]_q=\frac{[m]_q!}{[n]_q![m-n]_q!},\quad
  [n]_q!=[n]_q\cdot \ldots\cdot [1]_q,\quad
  [n]_q=\frac{q^n-q^{-n}}{q-q^{-1}}.
\end{equation*}

The comultiplication $\Delta$, the counit $\varepsilon$, and the
antipode $S$ are defined as follows:
\begin{equation*}
  \Delta(E_i)=E_i\otimes 1+K_i\otimes E_i,\, \,
  \Delta(F_i)=F_i\otimes K_i^{-1}+1\otimes F_i, \, \,
  \Delta(K_i)=K_i\otimes K_i,
\end{equation*}
\begin{equation*}
  S(E_i)=-K_i^{-1}E_i,\quad
  S(F_i)=-F_iK_i,\quad
  S(K_i)=K_i^{-1},
\end{equation*}
$$
\varepsilon(E_i)=\varepsilon(F_i)=0,\quad \varepsilon(K_i)=1.
$$

A representation $\rho: U_q \mathfrak{g} \rightarrow {\rm End} V$ is
called {\it weight} (and $V$ is called a weight module,
respectively), if $V$ admits a decomposition into a sum of weight
subspaces
\begin{align*}
V = \bigoplus_{\lambda} V_{\lambda}, \qquad V_{\lambda}=\{v \in V
|\text{ } \rho (K_j^{\pm 1})v = q_j^{\pm \lambda_j}v, j=1,...,l \},
\end{align*}
The subspace $V_{\lambda}$ is called a weight subspace of weight
$\lambda$.

It is convenient to define linear operators $H_i$ in a weight module
$V$ by
$$
H_iv=jv \qquad \mathrm{iff} \qquad K_iv=q_i^{j}v, \qquad v \in V.
$$

Let $U_q \mathfrak{k} \subset U_q \mathfrak{g}$ be a Hopf subalgebra
generated by $E_i,F_i, i=1,...,l, i \neq l_0$ and $K^{\pm 1}_j,
j=1,...,l$. A finitely generated weight $U_q \mathfrak{g}$-module $V$
is called a {\it quantum Harish-Chandra module} if $V$ is a sum of
finite dimensional simple $U_q \mathfrak{k}$-modules and $\dim
\rm{Hom}_{U_q \mathfrak{k}}(W,V)<\infty$ for every finite dimensional
simple $U_q \mathfrak{k}$-module $W$.

We restrict our consideration to quantum Harish-Chandra modules only.

Let $\lambda \in P_+$. $L(\lambda)$ denotes the simple $U_q
\mathfrak{g}$-module with a single generator $v(\lambda)$ and the
defining relations (see \cite{Jantzen})
$$
K_j^{\pm}v(\lambda)=q^{\pm \lambda_j}v(\lambda), \quad
E_jv(\lambda)=0, \quad F_j^{\lambda_j+1}v(\lambda)=0.
$$
Recall the notion of quantum analog of the algebra $\mathbb{C}[G]$ of
regular functions on $G$ \cite{Jantzen}. Denote by $\mathbb{C}[G]_q \subset
(U_q \mathfrak{g})^*$ the Hopf subalgebra of all matrix coefficients of
weight finite dimensional $U_q \mathfrak{g}$-representations. Denote by
$U_q^{op} \mathfrak{g}$ a Hopf algebra that differs from $U_q \mathfrak{g}$
by the opposite multiplication. Equip $\mathbb{C}[G]_q$ with a structure of
$U_q^{op} \mathfrak{g} \otimes U_q \mathfrak{g}$-module algebra in the
following way: $(\xi \otimes \eta)f{=}(L_{\rm reg}(\xi) \otimes R_{\rm
reg}(\eta))f$, where
$$
(R_{\rm reg}(\xi) f)(\eta)=f(\eta \xi), \qquad (L_{\rm
reg}(\xi)f)(\eta)=f(\xi\eta), \qquad \xi,\eta \in U_q \mathfrak{g}, f
\in \mathbb{C}[G]_q.
$$
The algebra $\mathbb{C}[G]_q$ is called the algebra of regular
functions on the quantum group $G$.

Introduce special notation for some elements of $\mathbb{C}[G]_q$
\cite{Soib}. Consider the finite dimensional simple $U_q
\mathfrak{g}$-module $L(\lambda)$ with highest weight $\lambda \in
P_+$. Equip it with an invariant scalar product $(\cdot,\cdot)$,
given by $(v(\lambda),v(\lambda))=1$ (as usual in the compact quantum
group theory \cite{Soib}). Choose an orthonormal basis of weight
vectors $\{v_{\mu,j}\} \in L(\lambda)_{\mu}$ for all weights $\mu$.
Let
$$
c^{\lambda}_{\mu,i,\mu',j}(\xi)= (\xi v_{\mu',j},v_{\mu,i}).
$$
We will omit the indices $i,j$, if it won't led to an ambiguity.
Introduce an auxiliary $U_q\mathfrak{g}$-module algebra
$\mathbb{C}[\widehat{X}]_q$.

Let $\lambda\in P_+$. For all $\lambda',\lambda''\in P_+$ $$ \dim
\operatorname{Hom}_{U_q \mathfrak{g}}
(L(\lambda'+\lambda''),L(\lambda')\otimes L(\lambda''))=1.
$$ Hence, the following $U_q\mathfrak{g}$-morphisms are
well-defined:
\begin{equation*}\label{m_lambda}
m_{\lambda',\lambda''}:L(\lambda')\otimes L(\lambda'')\to
L(\lambda'+\lambda''),\qquad
m_{\lambda',\lambda''}:v(\lambda')\otimes v(\lambda'')\mapsto
v(\lambda'+\lambda'').
\end{equation*}
Therefore, the vector space
\begin{equation*}\label{nondeg_cone}
\mathbb{C}[\widehat{X}]_q\overset{\mathrm{def}}{=}
\bigoplus_{\lambda\in P_+}L(\lambda)
\end{equation*}
is equipped with a $U_q\mathfrak{g}$-module algebra structure as follows:
$$
f'\cdot f''\overset{\mathrm{def}}{=}m_{\lambda',\lambda''}(f'\otimes
f''), \qquad \qquad f'\in L(\lambda'),\; f''\in L(\lambda'').
$$
This is a well-known quantum analog of the homogeneous coordinate
ring of the flag manifold $X=B \backslash G$.

The Peter-Weyl theorem claims that
$$
\mathbb{C}[G]_q=\bigoplus_{\lambda\in P_+}(L(\lambda)\otimes
L(\lambda)^*),
$$
Hence, there exists an embedding of $U_q\mathfrak{g}$-module algebras
\begin{equation*}\label{F_embed}
i:\mathbb{C}[\widehat{X}]_q\hookrightarrow\mathbb{C}[G]_q,\qquad
i:v(\lambda)\mapsto c_{\lambda,\lambda}^\lambda,\quad\lambda \in P_+,
\end{equation*}
where $c_{\lambda,\lambda}^\lambda$ are matrix coefficients of
representations introduced above, so we have

\begin{proposition}\label{int}
$\mathbb{C}[\widehat{X}]_q$ is an integral domain.
\end{proposition}

$\mathbb{C}[\widehat{X}]_q$ is a $P$-graded algebra by obvious
reasons:
$$
\mathbb{C}[\widehat{X}]_q=\bigoplus_{\lambda \in P_+}
\mathbb{C}[\widehat{X}]_{q,\lambda}, \qquad
\mathbb{C}[\widehat{X}]_{q,\lambda}=L(\lambda).
$$

\bigskip A simple weight finite dimensional $U_q\mathfrak{g}$-module
$L(\lambda)$ is called {\it spherical}, if it has a nonzero
$U_q\mathfrak{k}$-invariant vector \cite{Helg1}. This
$U_q\mathfrak{k}$-invariant vector is unique up to a nonzero
constant, $\dim L(\lambda)^{U_q \mathfrak{k}}=1$. Moreover, the set
$\Lambda_+\subset P_+$ of all highest weights of the spherical
$U_q\mathfrak{g}$-modules $L(\lambda)$ coincides with the similar set
in the classical case, so
$\Lambda_+=\bigoplus_{i=1}^r\mathbb{Z}_+\mu_i$, with
$\mu_1,\mu_2,\ldots,\mu_r$ being the fundamental spherical weights,
and $r$ being the real rank of bounded symmetric domain $\mathbb D$,
see \cite{Helg1}. Choose nonzero elements
$$\psi_i \in L(\mu_i)^{U_q \mathfrak{k}},\qquad i=1,2,\ldots,r.$$

\begin{proposition}\label{psi_commute}
$\psi_1,\psi_2,\ldots,\psi_r \in \mathbb{C}[\widehat{X}]_q$ pairwise
commute.
\end{proposition}

{\bf Proof.} $\psi_i\psi_j= \mathrm{const}(i,j)\psi_j\psi_i,$ with
$\mathrm{const}(i,j)\ne 0$, since $\dim (L(\mu_i+\mu_j)^{U_q
\mathfrak{k}})=1$. Prove that $\mathrm{const}(i,j)=1$. It follows
from Appendix, Lemma \ref{adj}, that we can equip
$\mathbb{C}[\widehat{X}]_q$ with an involution $\star$, and without
loss of generality we can assume $\psi_j=\psi_j^\star$ for all
$j=1,...,r$. Therefore, $\psi_i\psi_j=\pm\psi_j\psi_i$. The
involution $\star$ is a morphism of a continuous vector bundle
$\mathcal{E}_\lambda$ over $(0,1]$ with fibers $L(\lambda)_q$ defined
in Appendix. In the classical case $\psi_i\psi_j=\psi_j\psi_i$, so
the same holds in the quantum case. \hfill $\square$

\medskip
Consider a multiplicative subset
$$
\Psi=\{\psi_1^{j_1}\psi_2^{j_2}\cdots\psi_r^{j_r}|\:
j_1,j_2,\ldots,j_r\in\mathbb{Z}_+\}
$$
of the algebra $\mathbb{C}[\widehat{X}]_q$.

\begin{proposition}\label{psi_Ore}
$\Psi\subset\mathbb{C}[\widehat{X}]_q$ is an Ore set.
\end{proposition}

{\bf Proof.} Consider the decomposition of the
$U_q\mathfrak{g}$-module $L(\lambda)$ into a sum of its $U_q
\mathfrak{k}$-isotypic components $L(\lambda)=\bigoplus\limits_\mu
L(\lambda)_\mu$. Fix a subspace $L(\lambda)_\mu$ and a nonzero
$U_q\mathfrak{k}$-invariant vector $\psi\in L(\lambda')$. For all
$j\in\mathbb{N}$,
$$
L(\lambda+j\lambda')_\mu\supset L(\lambda+(j-1)\lambda')_\mu\,\psi.
$$
It follows from Proposition \ref{int} that
$L(\lambda+j\lambda')_\mu=L(\lambda+(j-1)\lambda')_\mu\,\psi$ for all
large enough $j\in\mathbb{N}$, since
\begin{equation}\label{mult_inequality}
\dim\operatorname{Hom}_{U_q\mathfrak{k}}(L(\mathfrak{k},\mu),L(\lambda))\le
\dim L(\mathfrak{k},\mu)
\end{equation}
for any simple weight finite dimensional $U_q\mathfrak{k}$-module
$L(\mathfrak{k},\mu)$ and any simple weight finite dimensional
$U_q\mathfrak{g}$-module $L(\lambda)$. The inequality
\eqref{mult_inequality} is well-known in the classical case
\cite[p.206]{Knapp86}, and the quantum case follows from the
classical case. Therefore,
$$\psi^j\,L(\lambda)_\mu\subset L(\lambda+(j-1)\lambda')_\mu\,\psi$$
for all large enough $j\in\mathbb{N}$, since $\psi^j\,
L(\lambda)_\mu\subset L(\lambda+j\lambda')_\mu$. Similarly,
$L(\lambda)_\mu\,\psi^j\subset\psi \,L(\lambda+(j-1)\lambda')_\mu$
for all large enough $j\in\mathbb{N}$. Hence, for all
$f\in\mathbb{C}[\widehat{X}]_q$, $\psi\in\Psi$
$$
\Psi f\cap\mathbb{C}[\widehat{X}]_q\,\psi\ne\varnothing,\qquad
f\Psi\cap\psi\,\mathbb{C}[\widehat{X}]_q\ne\varnothing,
$$
which is just the Ore condition for $\Psi$. \hfill $\square$

\bigskip Consider the localization $\mathbb{C}[\widehat{X}]_{q,\Psi}$
of the algebra $\mathbb{C}[\widehat{X}]_q$ with respect to the
multiplicative set $\Psi$.

The $P$-grading can be extended to
$\mathbb{C}[\widehat{X}]_{q,\Psi}$, since the elements of $\Psi$ are
homogeneous. The subalgebra
\begin{equation*}\label{nondeg_omega}
\mathbb{C}[\Omega]_q=\{f\in\mathbb{C}[\widehat{X}]_{q,\Psi}|\: \deg
f=0\}
\end{equation*}
is a quantum analog of the algebra $\mathbb{C}[\Omega]$ of regular
functions on the open $K$-orbit $\Omega \subset X=B\backslash G$.

\section{\boldmath $U_q\mathfrak{g}$-module algebra
$\mathbb{C}[\widehat{X}]_{q,\Psi}$} \label{U_cone}

In this section we equip $\mathbb{C}[ \widehat{X}]_{q,\Psi}$ with a $U_q
\mathfrak{g}$-module algebra structure. Start with some auxiliary facts.

Consider the vector space $L$ of sequences $\{x_n\}_{n\in
\mathbb{Z}_+}$
\begin{equation}\label{defL}
x_n=\sum_{i,j\in \mathbb{Z}_+}a_{ij}\lambda_i^nn^j,
\end{equation}
where only finite number of terms are nonzero, $a_{ij} \in
\mathbb{C}[\widehat{X}]_{q,\Psi}$, $\lambda_i$ are nonzero and
pairwise different. Sequences $x_n'$ and $x_n''$ are called {\it
asymptotically equal}, $x_n'\underset{\mathrm{as}}{=}x_n''$, if there
exists $N\in\mathbb{Z}_+$ such that for all $n \geq N$ one has
$x_n'=x_n''$.

\begin{lemma}\label{unique_alg}
Let $\{x_n\}_{n \in \mathbb{Z}_+} \in L$ and
$x_n\underset{\mathrm{as}}{=}0$. Then $x_n=0$ for all
$n\in\mathbb{Z}_+$.
\end{lemma}

{\bf Proof.} Let $x_n\underset{\mathrm{as}}{=}0$ and $x_n\not\equiv
0$. Without loss of generality, one assumes that $x_0\ne 0$, since
$L$ is invariant under $T$, with
$$
T:\{u_0,u_1,u_2,\ldots\}\mapsto\{u_1,u_2,\ldots\}, \qquad u_j \in
 \mathbb{C}[\widehat{X}]_{q,\Psi}.
$$
Let $\mathcal{M}$ be the smallest $T$-invariant subspace containing
all nonzero terms $\{a_{ij}\lambda_i^n n^j\}_{n\in \mathbb{Z}_+}$
from \eqref{defL}. Then $\dim\mathcal{M}<\infty$ and $ 0 \not \in
\mathrm{spec}(T|_\mathcal{M})$. However, $(T|_\mathcal{M})^N(x_n)=0$
for all large enough $N$, since $x_n\underset{\mathrm{as}}{=}0$. So,
$\mathrm{spec}(T|_\mathcal{M})=\{0\}$. This contradiction completes
the proof. \hfill $\square$

\begin{lemma} \label{Bernoulli}
Let $k \in \mathbb{Z}_+$ and $\lambda_0 \in \mathbb C$. The sequence
\begin{equation}\label{secuen}
\left\{\frac{d^k}{d\lambda^k}(1+\lambda+\lambda^2+\cdots+\lambda^n)_{
|_{\lambda=\lambda_0}}\right\}_{n \in \mathbb{Z}_+}
\end{equation}
belongs to $L$.
\end{lemma}

{\bf Proof.} Let $\lambda_0=1$. \eqref{secuen} is a sequence of
values of polynomials at $n \in \mathbb{Z}_+$, since
$$\sum_{j=0}^{n-1}j^{k-1}=\frac{B_k(n)-B_k(0)}{k},$$
with $B_k(z)$ being the Bernoulli polynomials \cite[chap.
3]{Prasolov}. For $\lambda_0 \neq 1$,
$$1+\lambda+\lambda^2+\cdots+\lambda^{n-1}=\frac{\lambda^n-1}{\lambda-1},$$
and the statement is now obvious.\hfill $\square$

\bigskip Let $\psi_0=\prod\limits_{j=1}^r\psi_j.$
Note that for any $f\in \mathbb{C}[\widehat{X}]_{q,\Psi}$, $\xi \in
U_q\mathfrak{g}$, the element $\xi(f\psi_0^n)$ is already defined and
belongs to $\mathbb{C}[\widehat{X}]_q$ for large enough $n \in
\mathbb Z_+$.

\begin{proposition}\label{U_actions}
For any $\xi\in U_q\mathfrak{g}$ there exists a unique linear
operator $R_\xi$ in $\mathbb{C}[\widehat{X}]_{q,\Psi}$ such that for
any $f\in\mathbb{C}[\widehat{X}]_{q,\Psi}$
\begin{itemize}
\item[1.] $R_\xi(f\psi_0^n)\underset{\mathrm{as}}{=}\xi(f\psi_0^n)$;

\item[2.] the sequence $x_n=R_\xi(f\psi_0^n)\cdot\psi_0^{-n}$ belongs
to $L$.
\end{itemize}
\end{proposition}

{\bf Proof.} The uniqueness follows from Lemma \ref{unique_alg}.
Prove the existence.

\medskip Consider a set of all $\xi\in U_q\mathfrak{g}$ such that there
exists $R_\xi$ with the required properties. The set forms a
subalgebra. Indeed,
$$\xi(f\psi_0^n)\underset{\mathrm{as}}{=}
R_\xi(f\psi_0^n)= \left(\sum_{i,j\in
 \mathbb{Z}_+}a_{ij}\lambda_i^nn^j\right)\psi_0^n,
$$
$$
\eta(a_{ij}\psi_0^n)\underset{\mathrm{as}}{=}
 R_\eta(a_{ij}\psi_0^n)=
\left(\sum_{k,l\in
 \mathbb{Z}_+}b_{ijkl}\mu_k^nn^l\right)\psi_0^n,
$$
so
$$
\eta\xi(f\psi_0^n)\underset{\mathrm{as}}{=} \left(\sum_{i,j,k,l\in
\mathbb{Z}_+}b_{ijkl}(\lambda_i\mu_k)^nn^{(j+l)}\right)\psi_0^n,
$$
and one can put $R_{\eta\xi}(f)=\sum_{i,k}\,b_{i0k0}$.
 From Lemma \ref{unique_alg},
$R_{\eta\xi}(f)$ is a linear operator in
$\mathbb{C}[\widehat{X}]_{q,\Psi}$. The sequence
$R_{\eta\xi}(f\psi_0^n)$ has all required properties.

\medskip
Now we have to construct linear operators
$$R_{K_i^{\pm 1}},\quad R_{E_i},\quad R_{K_iF_i},\qquad i=1,2,\ldots,l,
$$
that satisfy conditions of Proposition \ref{U_actions}. The case
$\xi=K_i^{\pm 1}$ is trivial, while the two others are very similar.

Prove the existence of $R_{E_i}$. Let $N\in \mathbb{Z}_+$ be the
smallest number such that $f\psi_0^N\in\mathbb{C}[\widehat{X}]_q$.
For any $n\in\mathbb{Z}_+$
\begin{gather*} E_i\left(f\psi_0^{n+N}\right)\psi_0^{-(n+N)}=
E_i\left(f\psi_0^N\right)\psi_0^{-N}+
\left(K_i\left(f\psi_0^N\right)\psi_0^{-N}\right)
\left(\psi_0^NE_i\left(\psi_0^n\right)\psi_0^{-(n+N)}\right).
\end{gather*}
It is enough to prove that $x_n=\psi_0^N\cdot
E_i(\psi_0^n)\cdot\psi_0^{-(n+N)}$ looks like \eqref{defL}, since
then one can put
$$
E_i(f\psi_0^n)\underset{\mathrm{as}}{=}
 \left(\sum_{s,t\in \mathbb{Z}_+}a_{st}\lambda_s^nn^t\right)\psi_0^n,
$$
and $R_{E_i}(f)=\sum_{s \in \mathbb Z_+}a_{s0}$. The linearity of
$R_{E_i}$ and all required properties reduce to the special case
$f\in\mathbb{C}[\widehat{X}]_q$ by Lemma \ref{unique_alg}, as before.

Turn to the proof. Let $A$ be the linear operator in
$\mathbb{C}[\widehat{X}]_{q,\Psi}$, defined by
$$
Af=\psi_0 f\psi_0^{-1}, \qquad f \in
\mathbb{C}[\widehat{X}]_{q,\Psi}.
$$
The key observation is that the linear span of
$\{A^k(E_i\psi_0\}_{k\in\mathbb{Z}_+}$ is finite dimensional. Indeed,
all vectors $A^k(E_i\psi_0)$ belong to the same homogeneous component
of the $P$-graded vector space $\mathbb{C}[\widehat{X}]_{q,\Psi}$,
and, moreover, to the same its $U_q\mathfrak{k}$-isotypic component
that corresponds to the highest weight $\mu=\sum\limits_{j=1}^r
j\mu_j+\alpha_i$. But the dimension of intersection $\mathcal{E}$ of
these components is at most than $\dim L(\mathfrak{k},\mu)^2$ by
\eqref{mult_inequality}.

Now one can use the equality
$$
x_n=\psi_0^N\left(\sum_{j=0}^{n-1}\psi_0^j(E_i\psi_0)\psi_0^{n-1-j}\right)
\psi_0^{-(n+N)}= \psi_0^N
\left(\sum_{j=0}^{n-1}(A|_{\mathcal{E}})^j\,(E_i\psi_0)\right)
\psi_0^{-(N-1)}
$$
and Lemma \ref{Bernoulli} to compute the function
$1+z+z^2+\cdots+z^{n-1}$ of the finite dimensional linear operator
$A|_{\mathcal{E}}$ reduced to the Jordan canonical form \cite{Mal}.
\hfill $\square$

\medskip
\begin{corollary}
$R_{\xi}$ provides a $U_q\mathfrak{g}$-module structure in
$\mathbb{C}[\widehat{X}]_{q,\Psi}$ which extends the
$U_q\mathfrak{g}$-module structure from $\mathbb{C}[\widehat{X}]_q$
to $\mathbb{C}[\widehat{X}]_{q,\Psi}$.
 \end{corollary}

\medskip

Now prove that $\mathbb{C}[\widehat{X}]_{q,\Psi}$ is a
$U_q\mathfrak{g}$-module algebra. Let $\widetilde{L}$ be the vector
space of functions on $\mathbb{Z}^2$ that takes values in
$\mathbb{C}[\widehat{X}]_{q,\Psi}$, such that
$$
x(m,n)= \sum_{i',i'',j',j''\in
\mathbb{Z}_+}a_{i'i''j'j''}\lambda_{i'}^n\,n^{j'}\,\mu_{i''}^m\,m^{j''},
$$
where the sum is finite, $a_{i'i''j'j''} \in
\mathbb{C}[\widehat{X}]_{q,\Psi}$, and $\lambda_{i'}, \mu_{i''} \in
\mathbb{C}$. Suppose that all $\lambda_{i'}$, $\mu_{i''}$ are
pairwise different and nonzero. One can easy expand Lemma
\ref{unique_alg} to functions $x(m,n) \in \widetilde{L}$. Namely, if
$x(m,n)=0$ for any $m\ge M$, $n\ge N$, then $x(m,n)=0$ for all
$m,n\in\mathbb{Z}$.

It is important to note that both in the statement and in the proof
of Proposition \ref{U_actions} one can replace the conditions on
$R_\xi$ by the following conditions:
\begin{itemize}
\item[1.] $R_\xi(\psi_0^mf\psi_0^n)=\xi(\psi_0^mf\psi_0^n)$, if
$m,n\in\mathbb{Z}_+$ and $m+n$ is large enough;

\item[2.] the function
$x(m,n)=\psi_0^{-m}R_\xi(\psi_0^mf\psi_0^n)\psi_0^{-n}$ belongs to
$\widetilde{L}$.
\end{itemize}
One gets the same representation $R_\xi$ of $U_q\mathfrak{g}$ in the
vector space $\mathbb{C}[\widehat{X}]_{q,\Psi}$.

\bigskip
\begin{proposition}\label{they_respect}
Let $\xi\in U_q\mathfrak{g}$ and
$\triangle\xi=\sum\limits_i\xi_i'\otimes\xi_i''$. Then for any
$f_1,f_2\in\mathbb{C}[\widehat{X}]_{q,\Psi}$ one has
$$\xi(f_1f_2)=\sum\limits_i(\xi_i'f_1)(\xi_i''f_2).$$
\end{proposition}

{\bf Proof.} Let
$$
x'(m,n)=\psi_0^{-m}\xi(\psi_0^mf_1f_2\psi_0^n)\psi_0^{-n},
$$
and
$$
x''(m,n)=
\sum_i\psi_0^{-m}\xi_i'(\psi_0^mf_1)\xi_i''(f_2\psi_0^n)\psi_0^{-n}.
$$

They are equal for large enough $m,n$, since
$\psi_0^mf_1,f_2\psi_0^n\in\mathbb{C}[\widehat{X}]_q$, and
$\mathbb{C}[\widehat{X}]_q$ is a $U_q\mathfrak{g}$-module algebra.
Also, both belong to $\widetilde{L}$. Therefore, $x'(m,n)=x''(m,n)$
for any $m,n\in\mathbb{Z.}$ Now put $m=n=0$ in the last equality.
\hfill $\square$

\begin{remark}
It is worth to note that Lunts and Rosenberg described another
approach to extension the $U_q\mathfrak{g}$-module algebra structure
in \cite{LuntsRosLoc},\cite{LuntsRosDiff2}. Their approach is more
general but more intricate.
\end{remark}

\section{Degenerate and non-degenerate spherical principal series}
\label{principal_constr} For simplicity and more clear presentation,
start with certain degenerate spherical principal series. The same
approach is used in producing non-degenerate spherical principal
series.

Fix $k \in \{1,2,\cdots,r\}$. Consider
$$\mathbb{C}[\widehat{X}_k]_q\overset{\mathrm{def}}{=}
\bigoplus_{j\in\mathbb{Z}_+}L(j\mu_k).
$$
As in the previous section, one equips $\mathbb{C}[\widehat{X}_k]_q$
with a $U_q\mathfrak{g}$-module algebra structure.
$\mathbb{C}[\widehat{X}_k]_q$ naturally embeds in the
$U_q\mathfrak{g}$-module algebra $\mathbb{C}[\widehat{X}]_q$ and has
a $\mathbb Z_+$-grading:
$$
\deg f =j, \qquad \text{iff} \qquad f \in L(j\mu_k).
$$
It follows from Proposition \ref{psi_Ore} that
\hbox{$\Psi_k=\psi_k^{\mathbb{Z}_+}$} is an Ore set, and the
localization $\mathbb{C}[\widehat{X}_k]_{q,\Psi_k}$ is a $\mathbb
Z$-graded $U_q \mathfrak{g}$-module algebra. Consider the subalgebra
$$\mathbb{C}[\Omega_k]_q=\{f\in\mathbb{C}[\widehat{X}_k]_{q,\Psi_k}|\:
\deg(f)=0\}.
$$
Evidentally, $\mathbb{C}[\Omega_k]_q \subset \mathbb{C}[\Omega]_q$ is
a $U_q\mathfrak{g}$-module algebra.

For $u \in\mathbb{Z}$ denote
\begin{equation}\label{pi_u}
\pi_{k,u}(\xi)f \stackrel{\rm def}{=}
\xi\left(f\psi_k^u\right)\cdot\psi_k^{-u}, \qquad\xi\in
U_q\mathfrak{g},\quad f\in\mathbb{C}[\Omega_k]_q.
\end{equation}
The representations $\pi_{k,u}$ are representations of degenerate
spherical principal series. Now we are going to introduce
$\pi_{k,u}$ for arbitrary $u \in \mathbb C$. We need some auxiliary
constructions.

Let
$$
\mathbb{C}[\Omega_k]_q=\bigoplus_{\lambda\in
P_+^\mathcal{S}}\mathbb{C}[\Omega_k]_{q,\lambda}
$$
be a decomposition of $\mathbb{C}[\Omega_k]_q$ into a sum of its
$U_q\mathfrak{k}$-isotypic components. $P_+^\mathcal{S}$ denotes the
set of all integral dominant weights of $\mathfrak{k}$. By
considerations from Appendix, $\mathbb{C}[\Omega_k]_{q,\lambda}$ are
fibers of a continuous vector bundle $\mathcal{F}_\lambda$ over
$(0,1]$ that is analytic on $(0,1)$. We identify morphisms of such
vector bundles over $(0,1]$ with the corresponding continuous in
$(0,1]$ and analytic in $(0,1)$ "operator valued functions".

It is easy to prove that the operator valued function
$$
A_{k,\lambda}(q):\mathbb{C}[\Omega_k]_{q,\lambda}\to
\mathbb{C}[\Omega_k]_{q,\lambda},\qquad
A_{k,\lambda}(q):f\mapsto\psi_kf\psi_k^{-1}
$$
is well-defined, invertible, continuous in $(0,1]$ and analytic in $(0,1)$.

\begin{lemma}\label{real_spectr}
All eigenvalues of $A_{k,\lambda}(q)$ are positive, and, moreover, rational
powers of $q$.
\end{lemma}

{\bf Proof.} Let $a$ be an eigenvalue of $A_{k,\lambda}(q)$, thus
there exists a nonzero $f \in L(j\mu_k) \subset
\mathbb{C}[\widehat{X}_k]_q$, such that $ \psi_k f=a\, f \psi_k$.

It can be shown easily that $aq^{-j(\mu_k,\mu_k)}$ is an eigenvalue
of the linear operator $R_{L(j\mu_k) L(\mu_k)}$ corresponded to the
universal $R$-matrix of $U_q \mathfrak{g}$. Here $(\cdot,\cdot)$ is
fixed by $(\alpha_i,\alpha_j)=d_ia_{ij}$.

It remains to prove that eigenvalues of $R_{L(\lambda),L(\lambda')}$ are
rational powers of $q$ for all $\lambda, \lambda' \in P_+$. There exists a
suitable basis of tensor products of weight vectors such that the matrix of
$R_{L(\lambda),L(\lambda')}$ is upper-triangular and its diagonal elements
belong to the set $ \{ q^{-(\mu',\mu'')} \,|\; \mu',\mu'' \in P \}$, hence,
they are rational powers of $q$. \hfill $\square$

\medskip
$\pi_{k,u}(K_i^{\pm 1}), \pi_{k,u}(E_i), \pi_{k,u}(F_i)$ are defined
for $u \in \mathbb Z$. We are going to extend these operator valued
functions to the complex plane. Evidently,
\begin{gather*}
\pi_{k,u}(K_i^{\pm 1})f=K_i^{\pm 1}(f \psi_k^u) \psi_k^{-u}=K_i^{\pm 1} f
K_i^{\pm 1} (\psi_k^u) \psi_k^{-u},
\\ \pi_{k,u}(E_i)f= E_i f+K_if E_i(\psi_k^u) \psi_k^{-u},
\\ \pi_{k,u}(K_iF_i)f = K_iF_if + K_if K_iF_i (\psi_k^u) \psi_k^{-u},
\qquad f \in \mathbb{C}[\Omega_k]_q.
\end{gather*}

Denote by $P_\lambda$ a projection in $\mathbb{C}[\Omega_k]_q$ with $
\rm{Im} \,P_\lambda=\mathbb{C}[\Omega_k]_{q,\lambda}$, $\rm{Ker}
\,P_\lambda= \bigoplus_{\lambda' \neq \lambda}
\mathbb{C}[\Omega_k]_{q,\lambda'}. $ In the sequel we deal with
operator valued functions
\begin{equation*}
P_{\lambda_2}\pi_{k,u}(K_i^{\pm
1})|_{\mathbb{C}[\Omega_k]_{q,\lambda_1}},\qquad
P_{\lambda_2}\pi_{k,u}(E_i)|_{\mathbb{C}[\Omega_k]_{q,\lambda_1}},\qquad
P_{\lambda_2}\pi_{k,u}(K_iF_i)|_{\mathbb{C}[\Omega_k]_{q,\lambda_1}}.
\end{equation*}

Consider a sequence $\{a_u\}_{u \in \mathbb Z_+}$
$$
a_u=E_i\psi_k^u=\left(\sum_{j=0}^{u-1}A^j(E_i\psi_k)\right) \psi_k^{u-1},
\qquad u \in \mathbb{Z}_+,
$$
where $A_k: \mathbb{C}[\Omega_k]_q \rightarrow
\mathbb{C}[\Omega_k]_q$ and
$A_k|_{\mathbb{C}[\Omega_k]_{q,\lambda}}=A_{k, \lambda}(q).$

Similarly to the arguments in the previous section, show that $\{E_i
\psi_k^u\}_u \in L$.

It is clear that $E_i\psi_k$ belongs to $V=\oplus_{\lambda \in
\mathcal{M}} \mathbb{C}[\Omega_k]_{q,\lambda}$, with $\mathcal{M}
\subset P_+^\mathcal{S}$ being a finite set. Consider the restriction
$A_k$ to $V$, $A_kV \subset V$. Using the Jordan canonical form of
$A_k|_V$ and the equality
\begin{equation*}
E_i\psi_k^u=\left(\sum_{j=0}^{u-1}A_k^j(E_i\psi_k)\right)
\psi_k^{u-1}, \qquad u \in \mathbb{Z}_+,
\end{equation*}
one has $\{E_i \psi_k^u\}_u \in L$.

So one can extend the operator valued function $\pi_{k,u}(E_i)$ to
the complex plane, since the eigenvalues of $A_{k,\lambda}(q)$ are
positive (Lemma \ref{real_spectr}). Similarly, one can extend the
operator valued function $\pi_{k,u}(K_iF_i)$. The extensions of the
operator valued functions $\pi_{k,u}(K^{\pm 1}_i)$, $i=1,2,\cdots,l$,
exist by obvious reasons. At last,
$$\pi_{k,u}(F_i)=\pi_{k,u}(K_i^{-1})\,\pi_{k,u}(K_iF_i), \qquad i=1,2,\cdots,l.$$

Now we have to check whether the map
$$
E_i \mapsto \pi_{k,u}(E_i), \quad F_i \mapsto \pi_{k,u}(F_i), \quad
K_i \mapsto \pi_{k,u}(K_i),
$$ can be extended to an algebra homomorphism.

Introduce an auxiliary algebra of analytic functions $F(u;q)$ on
$\mathbb{C}\times(0,1)$ that take values in the space of linear operators
in $\mathbb{C}[\Omega_k]_q=\bigoplus_{\lambda\in
P_+^\mathcal{S}}\mathbb{C}[\Omega_k]_{q,\lambda}$. In other words, we
assume the analyticity of all operator valued functions
$P_{\lambda_2}F(u;q)_{|\mathbb{C}[\Omega_k]_{q,\lambda_1}}$, where
$\lambda_1, \lambda_2 \in P_+^\mathcal{S}$.

The vector bundle $\mathcal{F}_\lambda$ is equipped with a Hermitian
metric, so the operator norm $\|F(u;q)\|$ is well-defined. Consider a
subalgebra of operator valued functions satisfying the following
condition:
\begin{equation}\label{ineq_unique}
\|F(u;q)\|\le a_{_F}(q)\exp(b_{_F}(q)|u|),
\end{equation}
for some $a_{_F}(q)>0$ and $b_{_F}(q)>0$ such that $\lim\limits_{q\to
1}b_{_F}(q)=0$. Note that the subalgebra does not depend on the
choice of metrics.

The operator valued functions
$$
\pi_{k,u}(E_i), \; \pi_{k,u}(F_i),\; \pi_{k,u}(K_i^{\pm 1}), \qquad
i=1,2,\cdots,l,
$$ are analytic and satisfy \eqref{ineq_unique}. Prove it. Consider $A_k|_V$, where $V$ is the
finite dimensional $U_q \mathfrak k$-invariant subspace and $E_i
\psi_k \in V$. One has
$$
\|(A_k|_V)^u(E_i\psi_k)\|=\|\sum_{s,t \in \mathbb Z_+} a_{st}
\lambda_s^u u^t\|.
$$
The number of terms in the r.h.s. of expression is at most $(\dim
V)^2$ (obviously, $0 \leq s,t \leq \dim V-1$). Hence,
$$
\|(A_k|_V)^u(E_i\psi_k)\|=\|\sum_{s,t \in \mathbb Z_+} a_{st}
\lambda_s^u u^t\| \leq (\dim V)^2 \max \|a_{st}\| \max|\lambda_s|^u
u^t.
$$

\begin{proposition}(\cite{Evgrafov})\label{unique_f}
Let $f(z)$ be continuous in $\{z \in \mathbb{C}| \operatorname{Re} z
\geq 0\}$ and holomorphic in $\{z \in \mathbb{C}| \operatorname{Re} z
> 0\}$. Assume also that
\begin{itemize}
\item[1.]
$|f(iy)|\le\mathrm{const}\cdot\exp\{(\pi-\varepsilon)|y|\}$,\ \ \
$y\in\mathbb{R}$,
 \item[2.] $|f(z)|\le M\exp(a|z|)$, \ \ \ $\operatorname{Re}z \ge 0$
\end{itemize}
for some real $M,a$ and $\varepsilon > 0$.

If $f(n)=0$ for all $n\in\mathbb{N}$, then $f(z)\equiv 0$.
\end{proposition}

\begin{corollary}
The extension of $\pi_{k,u}(E_i), \pi_{k,u}(F_i), \pi_{k,u}(K_i)$ is
unique.
\end{corollary}
Using \eqref{ineq_unique}, one can prove easily that the
Drinfeld-Jimbo relations hold for these operator valued functions. So
$\pi_{k,u}$ is a representation. It is a $q$-analog of a
representation of degenerate spherical principal series of
Harish-Chandra modules.

Now turn to the non-degenerate spherical principal series. For
$\mathbf{u}=(u_1, u_2,\cdots, u_r) \in \mathbb{Z}^r$ define (cf.
\eqref{pi_u})
\begin{equation*}
\pi_{\mathbf u}(\xi)f \stackrel{\rm def}{=} \xi \left( f\prod_{j=1}^r
\psi_j^{u_j}\right)\cdot\prod_{j=1}^r\psi_j^{-u_j}, \qquad \xi\in
U_q\mathfrak{g},\quad f\in\mathbb{C}[\Omega]_q.
\end{equation*}

We describe the extension of $\pi_{\mathbf u}$ to $\mathbb{C}^r$. As
before, consider the decomposition of $\mathbb{C}[\Omega]_q$ into a
sum of its $U_q\mathfrak{k}$-isotypic components
$$
\mathbb{C}[\Omega]_q=\bigoplus_{\lambda\in
P_+^\mathcal{S}}\mathbb{C}[\Omega]_{q,\lambda}
$$
and operator valued functions
$$
\mathcal{A}_{j,\lambda}(q):\mathbb{C}[\Omega]_{q,\lambda}\to
\mathbb{C}[\Omega]_{q,\lambda},\qquad
\mathcal{A}_{j,\lambda}(q):f\mapsto\psi_jf\psi_j^{-1}.
$$
The construction of
\begin{equation}\label{pi_final}
\pi_\mathbf{u}(E_i),\;\; \pi_\mathbf{u}(F_i),\;\;
\pi_\mathbf{u}(K_i^{\pm 1}), \qquad i=1,2,\cdots,l,
\end{equation}
essentially reduces to analytic continuation of the vector-valued
functions
$$
 \sum_{j=1}^{u_k-1}\mathcal{A}_{k,\lambda}^j
 (E_i\psi_k), \qquad k=1,2,\cdots,r.
$$

The Drinfeld-Jimbo relations for the operator valued functions
\eqref{pi_final} can be proved in the same way. Namely, consider an
algebra of analytic operator valued functions
$F(u_1,u_2,\cdots,u_r;q)$ such that
$$\|F(u_1,u_2,\cdots,u_r;q)\|\le a_{_F}(q)\exp(b_{_F}(q)\sum_{k=1}^r|u_k|),$$
(cf. \eqref{ineq_unique}). Now one can prove the uniqueness of the
interpolation of $\pi_{\mathbf u}$ in this subalgebra.

\section{Appendix}

We present some auxiliary statements on certain vector bundles over
$(0,1]$. Start with well-known facts on Verma modules. Let $U_q
\mathfrak{b}^+$ be a Hopf subalgebra generated by $E_i, K_i^{\pm 1}$.
Let $\lambda \in P_+$, and $\mathbb C_\lambda$ a one dimensional
$U_q\mathfrak{b}^+$-module defined by its generator $1_\lambda$ and
the relations
\begin{equation*}
K_i^{\pm 1} 1_\lambda = q_i^{\pm\lambda_i} 1_\lambda, \qquad E_i
1_\lambda = 0,\qquad i=1,2,\ldots,l.
\end{equation*}
As usual, a Verma module over $U_q \mathfrak{g}$ can be defined in a
such way:
$$
M(\lambda)_q\stackrel{\rm
def}{=}U_q\mathfrak{g}\otimes_{U_q\mathfrak{b}^+}\mathbb{C}_\lambda.
$$

Fix $v_\lambda=1\otimes 1_\lambda$. It is well-known that $v_\lambda$
is a generator, and $M(\lambda)_q$ can be defined by the relations:
$$
E_iv_\lambda=0,\qquad K_i^{\pm
1}v_\lambda=q_i^{\pm\lambda_i}v_\lambda,\qquad i=1,2,\ldots,l.
$$

Recall that the Weyl group $W$ acts on the root system R of Lie
algebra $\mathfrak{g}$ and is generated by simple reflections
$s_i(\alpha_j)=\alpha_j-a_{ij}\alpha_i$. Fix the reduced expression
of the longest element $w_0=s_{i_1}\cdot s_{i_2}\cdot \ldots \cdot
s_{i_M} \in W$. One can associate to it a total order on the set of
positive roots of $\mathfrak{g}$, and then a basis in the vector
space $U_q \mathfrak{g}$. We use the following total order on the set
of positive roots:
$$
\beta_1=\alpha_1,\quad \beta_2=s_{i_1}(\alpha_{i_2}),\quad
\beta_3=s_{i_1}s_{i_2}(\alpha_{i_3}),\quad \ldots
\quad\beta_M=s_{i_1}\ldots s_{i_{M-1}}(\alpha_{i_M}).
$$

Following Lusztig \cite{Lust1,Rosso,DaDeCon}, introduce elements
$E_{\beta_s}, F_{\beta_s} \in U_q \mathfrak{g}$ for $s=1,...,M-1$. As
a direct consequence of definitions $E_{\beta_s}$, (resp.
$F_{\beta_s}$) is a linear combination of $E_{j_1}^{m_1}\cdot ...
\cdot E_{j_l}^{m_l}$ (resp. $F_{i_1}^{n_1} \cdot...\cdot
F_{i_k}^{n_k}$) with the coefficients in the expansion being rational
functions of $q$ without poles in $(0,1]$.

\medskip

\begin{proposition}\label{la}
The set $\{F_{\beta_M}^{j_M}\cdot F_{\beta_{M-1}}^{j_{M-1}}\cdot
\ldots \cdot F_{\beta_1}^{j_1}\cdot K_1^{i_1}\cdot K_2^{i_2}\cdot
\ldots \cdot K_{l}^{i_l}\cdot E_{\beta_1}^{j_1}\cdot
E_{\beta_2}^{j_2}\cdot \ldots \cdot E_{\beta_M}^{j_M}|
\\ k_1,k_2,\ldots,k_M,\,j_1,j_2,\ldots,j_M \in \mathbb
{Z}_+\,, i_1,i_2,\ldots,i_l \in \mathbb{Z}\}$ is a basis in the
vector space $U_q \mathfrak{g}$.
\end{proposition}

Hence, the weight vectors
\begin{equation}\label{PBV-Verma}
v_J(\lambda)=F_{\beta_M}^{j_M}F_{\beta_{M-1}}^{j_{M-1}}\ldots
F_{\beta_1}^{j_1}v_\lambda,\qquad j_1,j_2,\ldots,j_M \in
\mathbb{Z}_+,
\end{equation}
form a basis of $M(\lambda)_q$.

Equip $U_q\mathfrak{g}$ with a $*$-Hopf algebra structure as follows:
\begin{equation*}
(K_j^{\pm 1})^\star=K_j^{\pm 1},\quad E_j^\star=K_jF_j,\quad
F_j^\star=E_jK_j^{-1},\quad j=1,2,\ldots,l.
\end{equation*}

\begin{lemma}
There exists a unique Hermitian form in $M(\lambda)_q$ such that
\begin{itemize}
\item $(\xi v',v'')=(v',\xi^\star v''),\qquad v',v''\in
M(\lambda)_q,\quad\xi\in U_q\mathfrak{g}$, \item
$(v_\lambda,v_\lambda)=1$.
\end{itemize}
\end{lemma}
The kernel $K(\lambda)_q$ of the form $(\cdot, \cdot)$ is the largest
proper submodule of $M(\lambda)_q$.

In this section we write $L(\lambda)_q$ instead of $L(\lambda)$ to
make the dependence on $q$ explicit.
\begin{proposition}
1. $L(\lambda)_q \simeq M(\lambda)_q/K(\lambda)_q$.

2. The form $(\cdot,\cdot)$ is non-degenerate in $L(\lambda)_q$.
\end{proposition}

Introduce a morphism $p_\lambda:M(\lambda)_q \rightarrow
L(\lambda)_q$, $v_\lambda \mapsto v(\lambda)$.

\medskip

Here is our first statement on special vector bundles over $(0,1]$.
Let $\lambda \in P_+$. There exists a continuous vector bundle
$\mathcal{E}_\lambda$ over $(0,1]$ with fibers isomorphic to
$L(\lambda)_q$, that is analytic in $(0,1)$ and $E_i, F_i, H_i$,
$i=1,2,\cdots,l$ act by endomorphisms in
$\mathcal{E}_\lambda$.\footnote{I.e. $E_i, F_i, H_i$,
$i=1,2,\cdots,l$ corresponds to continuous in $(0,1]$ and analytic in
$(0,1)$ operator valued functions.}

Describe the construction of $\mathcal{E}_\lambda$. Recall that to
any reduced expression of $w_0$ we assign the basis of
$M(\lambda)_q$, see \eqref{PBV-Verma}.

Fix $q_0\in (0,1]$. Choose a subset $\{v_j\}_{j=1,\ldots,\dim L(\lambda)}$
of one of the mentioned bases in a such way that the matrix
$((v_i,v_j))_{i,j=1,\ldots,\dim L(\lambda)}$ is non-degenerate. It is
non-degenerate for all $q$ that are close enough to $q_0$. Hence, in a
neighborhood of $q_0$, $\{p_\lambda(v_j)\}_{j=1,...,\dim L(\lambda)}$ is a
basis, since $\dim L(\lambda)_q$ does not depend on $q\in (0,1]$. One gets
a trivial vector bundle with the required properties over a neighborhood of
$q_0$.

Elements of the matrix are continuous in $(0,1]$ and analytic in
$(0,1)$ functions. Therefore, the matrices of $E_i, F_i, H_i$ in the
basis $\{p_\lambda(v_j)\}_{j=1,...,\dim L(\lambda)}$ are continuous
in $(0,1]$ and analytic in $(0,1)$. Indeed, any function
$$
(p_\lambda(E_{j_1}^{m_1}\cdot ... \cdot E_{j_l}^{m_l} F_{i_1}^{n_1}
\cdot...\cdot F_{i_k}^{n_k}v_\lambda), p_\lambda(v_\lambda))
$$
is continuous in $(0,1]$ and analytic in $(0,1)$, since to calculate
the value one should just use the commutation relations, which are
''well-dependent'' on $q$. Therefore, the functions
$$
(E_ip_\lambda(F_{i_1}^{n_1} \cdot...\cdot F_{i_k}^{n_k}v_\lambda),
p_\lambda(F_{j_1}^{m_1} \cdot...\cdot F_{j_k}^{m_k}v_\lambda))
$$
and
$$
(E_ip_\lambda(F_{\beta_1}^{n_1} \cdot...\cdot
F_{\beta_M}^{n_M}v_\lambda), p_\lambda(F_{\beta_1}^{m_1}
\cdot...\cdot F_{\beta_M}^{m_M}v_\lambda))
$$
are ''well-dependent'' on $q$. Hence, the matrix elements of the
operator $E_i$ are ''well-dependent'' on $q$, since the matrix is
non-degenerate in a neighborhood of $q_0$. The same holds for
$F_i,H_i$, and the transition matrices defined on intersection of the
neighborhoods. Finally, the vector bundle over $(0,1]$ which we
obtain in this way does not depend on the choices made above. So the
vector bundle $\mathcal{E}_\lambda$ is constructed.

\medskip

Now proceed to a construction of a subbundle of $\mathcal{E}_\lambda$
corresponding to a fixed $U_q \mathfrak{k}$-type $\mu$. Consider the
decomposition of $L(\lambda)_q=\bigoplus_\mu L(\lambda)_{q,\mu}$ into
a sum of its $U_q \mathfrak{k}$-isotypic components. Then one has the
following statement on special vector bundles. For any $\lambda \in
P_+$ and $\mu \in P_+^\mathcal{S}$, $L(\lambda)_{q,\mu}$ is a fiber
of a continuous vector bundle over $(0,1]$, analytic in $(0,1)$.

Fix $\mu \in P_+^\mathcal{S}$ and consider
$\mathcal{E}_\lambda^\mu=\{(f,q)|f \in L(\lambda)_{q,\mu}\}$. Prove
that it defines a subbundle of $\mathcal{E}_\lambda$. Indeed,
consider a fiber $L(\lambda)_q$ together with its decomposition
$L(\lambda)_q=\bigoplus_{\nu \in P_+^\mathcal{S}}L(\lambda)_{q,\nu}$.
Note that the sum consists of finite number of terms. Hence, there
exists $c_q \in Z(U_q \mathfrak{k})$, that is polynomial in $q$, and
$$
c_q|_{L(\lambda)_{q,\mu}}=1, \qquad c_q|_{L(\lambda)_{q,\nu}}=0,
\quad \nu \neq \mu,
$$
see \cite[p.125-126]{Jantzen}. Evidentally, $c_q$ defines a morphism
of vector bundles
$$\xymatrix{\mathcal{E}_\lambda \ar[dr]^j\ar[r]^{c_q} & \mathcal{E}_\lambda \ar[d]\\ &
(0,1]}$$ Also it is an orthogonal projection onto $L(\lambda)_{q,\mu}$ in
any fiber $L(\lambda)_q$. $\rm{rank} \,c_q$ is constant, so the image of
$c_q$ is a vector subbundle.

\medskip

Now we can construct the last required vector bundle.
$\mathbb{C}[\Omega_k]_{q,\lambda}$ are fibers of a continuous vector
bundle $\mathcal{F}_\lambda$ over $(0,1]$, analytic in $(0,1)$.

Consider a map
$$
\Psi_q: \mathbb{C}[\widehat{X}_k]_{q,\Psi_k} \rightarrow
\mathbb{C}[\widehat{X}_k]_{q,\Psi_k}, \qquad f \mapsto f\psi_k.
$$

It is easy to prove that $\Psi_q$ is an invertible, continuous
operator valued function, analytic in $(0,1)$. Using $\Psi_q$ one can
carry the vector bundle structure to $\mathcal{F}_\lambda$ from
$\{(f,q)|f \in \mathbb{C}[\Omega_k]_{q, \lambda}\psi_k^N , q \in
(0,1]\}$ for large enough $N$.

\medskip

%\begin{remark}
%Note that the matrix elements of all the operator valued functions
%belong to $\mathbb{Q}(q^{1/s})$ with $s=\rm{card}(P/Q)$, and $Q$
%being the root lattice.
%\end{remark}

The next considerations are related to self-adjointness. Consider an
auxiliary algebra
$$
\mathbb{C}[\widehat{X}^{\operatorname{spher}}]_q=\bigoplus_{\lambda
\in \Lambda_+} L(\lambda)_q.
$$
Equip $U_q \mathfrak{g}$ and
$\mathbb{C}[\widehat{X}^{\operatorname{spher}}]_q$ with a ''complex
conjugation''. Recall that $\dim L(\lambda)_q^{U_q \mathfrak{k}}=1$
for any $\lambda \in \Lambda_+$. Consider the antilinear involutive
automorphism $\bar{\cdot}$ of $U_q \mathfrak{g}$ defined by
$$ \bar{E}_i=E_i,\; \; \bar{F}_i=F_i,\; \;
\bar{K}_i^{\pm 1}=K_i^{\pm 1}, \qquad i=1,2,\cdots,l.$$ There exists
a unique antilinear involutive operator $\bar{\cdot}$
$$
L(\lambda)_q \rightarrow L(\lambda)_q, \qquad \xi v(\lambda) \mapsto
\bar{\xi}v(\lambda),\qquad \xi \in U_q\mathfrak{g}.
$$ (Indeed, the uniqueness is obvious, while the
existence follows from the definition of $L(\lambda)_q$.)

It is easy to see that $\overline{L(\lambda)_q^{U_q
\mathfrak{k}}}=L(\lambda)_q^{U_q \mathfrak{k}}$. Hence, there exists a
nonzero vector $w_\lambda \in L(\lambda)_q^{U_q \mathfrak{k}}$ such that
$w_\lambda=\bar{w_\lambda}$. Let $l(\lambda)=\mathbb R w_\lambda$.

\begin{lemma}\label{adj}
There exists a unique involution $\star$ of
$\mathbb{C}[\widehat{X}^{\operatorname{spher}}]_q$ such that
$\mathbb{C}[\widehat{X}^{\operatorname{spher}}]_q$ is a
$(U_q\mathfrak{g},\star)$-module algebra, and
$\star|_{l(\lambda)}=id$.
\end{lemma}
{\bf Proof.} Firstly prove that
\begin{equation}\label{self_duality}
L(\lambda)_q^*\approx L(\lambda)_q, \qquad \lambda \in \Lambda_+.
\end{equation}
Following \cite{Takeuchi}, introduce a system of strongly orthogonal roots
$\gamma_1 > \gamma_2 > \ldots > \gamma_r$ with $\gamma_1$ being the maximal
root. It is easy to prove that $-w_0(\gamma_j)=\gamma_j$ for
$j=1,2,\ldots,r,$ where $w_0\in W$ is the longest element. Hence,
\eqref{self_duality} follows from the fact that the fundamental spherical
weights belong to the linear span of $\gamma_1,\gamma_2,\ldots,\gamma_r$,
see \cite{Helg1}.

An involution $\star$ on $L(\lambda)_q$ such that
$$
(\xi f)^\star=(S(\xi))^\star f^\star,\qquad\xi\in
U_q\mathfrak{g},\quad f\in L(\lambda)_q,
$$
is unique up to $\pm 1$. The uniqueness of the involution now follows
from the fact that $w_\lambda=w_\lambda^\star$.

Turn to a proof of existence of $\star$. $U_q \mathfrak{g}$ is
equipped with the involution as follows:
$$
(K_j^{\pm 1})^\star=K_j^{\pm 1},\quad E_j^\star=K_jF_j,\quad
F_j^\star=E_jK_j^{-1},\quad j=1,2,\ldots,l.
$$
Consider the $*$-algebra $(\mathbb{C}[G]_q,\star)$ with the
involution $\star$ given by
$$
f^\star(\xi) \stackrel{\rm def}{=} \overline{f((S(\xi))^\star)},\quad
\xi\in U_q\mathfrak{g},\quad f \in \mathbb{C}[G]_q .
$$
Let
$$
F=\{f\in\mathbb{C}[G]_q\;|\;L_{\mathrm{reg}}(\xi)f=\varepsilon(\xi)f,
\quad\xi\in U_q\mathfrak{k}\}.
$$

It follows from Peter-Weyl expansion that $F
\approx\oplus_{\lambda\in\Lambda_+}L(\lambda)_q$ as a $U_q
\mathfrak{g}$-module. One can consider $\Lambda_+$ with a natural partial
order $\leq$, and $F$ can be equipped with a $(U_q
\mathfrak{g},\star)$-invariant filtration $F=\bigcup_{\lambda \in
\Lambda_+} F_\lambda,$ $F_\lambda=\bigoplus_{\mu \leq \lambda} L(\mu)_q$.

Consider zonal spherical functions $\rho_\lambda$ related to the
$U_q\mathfrak{g}$-modules $L(\lambda)_q$. Their images in the
associated $\mathbb{Z}^r$-graded algebra $\mathrm{Gr} F$ pairwise
commute. Indeed, the quasi-commutativity is evident, while the
commutativity follows from the self-adjointness of $\rho_\lambda$ and
the commutativity in the classical case. Now, it follows from
\cite[Cor. 2.6]{Vret} that $\mathrm{Gr} F$ is naturally isomorphic to
$\mathbb{C}[\widehat{X}^{\operatorname{spher}}]_q$. Note that
$l(\lambda)$ corresponds to $\mathbb{R} \rho_\lambda$. Carry the
involution $\star$ from $F$ to $\mathrm{Gr} F$ and
$\mathbb{C}[\widehat{X}^{\operatorname{spher}}]_q$. It can be
verified easily that $\star$ is a morphism of the vector bundle with
fibers $L(\lambda)_q$. \hfill $\square$
%\begin{remark}
%One can substitute the involution $\star$ with a certain involution
%$*$ such that $\mathbb{C}[\widehat{X}^{\operatorname{spher}}]_q$ is a
%$(U_q\mathfrak{g},*)$-module algebra, and $w_\lambda=w_{\lambda}^*$.
%In the previous considerations one should just substitute the
%involution $\star$ in $\mathbb{C}[G]_q$ with a suitable involution
%$*$, see \cite[chap. 2.3,3.3]{disser}. This way one gets a $(U_q
%\mathfrak{g},*)$-module algebra $\mathbb{C}[\Omega]_q$.
%We hope to
%obtain a maximal Furstenberg-Satake compactification with this $(U_q
%\mathfrak{g},*)$-module algebra, see \cite{Koranyi}.
%\end{remark}

\end{document}